\documentclass[12pt]{amsart}

\usepackage{amsmath}
\usepackage{amssymb}
\usepackage{amsthm}

\usepackage{fullpage}
\usepackage{tabularray}	
\usepackage{hyperref}

\usepackage{dynkin-diagrams}

\newtheorem{theorem}{Theorem}
\numberwithin{equation}{section}
\numberwithin{theorem}{section}

\newtheorem{lemma}[theorem]{Lemma}
\newtheorem{proposition}[theorem]{Proposition}

\theoremstyle{definition}
\newtheorem{definition}[theorem]{Definition}
\newtheorem{example}[theorem]{Example}

\theoremstyle{remark}
\newtheorem{remark}[theorem]{Remark}

\def\C{{\mathbb C}}
\def\R{{\mathbb R}}
\def\Z{{\mathbb Z}}

\def\cC{{\mathcal C}}
\def\cP{{\mathcal P}}

\def\cT{{\mathcal T}}

\def\g{{\mathfrak g}}
\def\h{{\mathfrak h}}
\def\b{{\mathfrak b}}
\def\i{{\mathfrak i}}
\def\p{{\mathfrak p}}
\def\t{{\mathfrak t}}
\def\m{{\mathfrak m}}
\def\r{{\mathfrak r}}
\def\z{{\mathfrak z}}
\def\n{{\mathfrak n}}
\def\e{{\mathfrak e}}

\def\G{{\mathfrak G}}
\def\B{{\mathfrak B}}
\def\P{{\mathfrak P}}
\def\0{_{\overline{0}}}
\def\1{_{\overline{1}}}

\def\sl{\mathfrak{sl}}
\def\so{\mathfrak{so}}
\def\psl{\mathfrak{psl}}
\def\osp{\mathfrak{osp}}

\def\Ab{\mathfrak{Ab}}
\def\Ad{\mathfrak{Ad}}

\DeclareMathOperator{\rA}{A}
\DeclareMathOperator{\rB}{B}
\DeclareMathOperator{\rC}{C}
\DeclareMathOperator{\rD}{D}
\DeclareMathOperator{\rF}{F}
\DeclareMathOperator{\rG}{G}

\DeclareMathOperator{\rank}{rank}
\newcommand{\admod}[1]{ad\hspace{2pt}{#1}-module}
\newcommand{\adsubmod}[1]{ad\hspace{2pt}{#1}-submodule}

\def\O{\varnothing}

\title{Ad-Nilpotent Ideals of Parabolic Subalgebras of Lie Superalgebras}
\author{Samuel Renforth}

\date{\today}

\thanks{The author receives funding through an EPSRC studentship.}

\begin{document}

\begin{abstract}
We prove that for every parabolic subalgebra $\P$ of a basic classical Lie superalgebra excluding $\psl(n|n)$, there exists a parabolic subalgebra $\p$ of a simple Lie algebra, such that the ad-nilpotent (respectively, abelian) ideals of $\p$ and $\P$ are in bijection. This is used to enumerate such ideals of Borel subalgebras of basic classical Lie superalgebras.
\end{abstract}

\maketitle

\section{Introduction}
In 1998, a remarkable result due to Peterson was recounted in~\cite{Ko1}, stating that a Borel subalgebra $\b$ of a finite-dimensional simple Lie algebra $\g$ over $\C$ has $2^{\rank(\g)}$ abelian ideals. This sparked interest in abelian ideals of $\b$, as well as in the broader class of ad-nilpotent ideals of $\b$, see for example~\cite{CP1,CP2,KOP,Pa1,Pa2,
Su} and the references therein.

A natural extension of Peterson's result is to enumerate ad-nilpotent and abelian ideals of parabolic subalgebras of $\g$. This was accomplished in 2008 by Righi in~\cite{Ri1,Ri2}.

In this paper we enumerate ad-nilpotent and abelian ideals of Borel subalgebras of basic classical Lie superalgebras. This is done by combining Righi's results with our main theorem.

\begin{theorem}\label{T:bij}
Let $\G$ be a basic classical Lie superalgebra excluding $\psl(n|n)$ and $\P$ a parabolic subalgebra of $\G$. There exists a parabolic subalgebra $\p$ of a simple Lie algebra $\g$, such that the ad-nilpotent (respectively, abelian) ideals of $\p$ and $\P$ are in bijection.
\end{theorem}
 
We explain why $\psl(n|n)$ is excluded in Remark~\ref{R:pslnn}. 
 
The number of ad-nilpotent and abelian ideals of a Borel subalgebra $\B$ of a basic classical Lie superalgebra $\G$ depends on $\B$, as the Weyl group of $\G$ does not in general act transitively on the sets of positive roots. 

Peterson's proof used properties of affine Weyl groups that do not appear to extend to Lie superalgebras. This is the reason for our approach.

Theorem~\ref{T:bij} is made possible due to Kostant's work on Kostant root systems of Lie algebras,~\cite{Ko2}, and Dimitrov and Fioresi's extension of these ideas to Lie superalgebras,~\cite{DF1}. It in particular builds on the proof of Theorem 1.4 of \textit{loc. cit.}, where it was shown that Lie superalgebra root systems are equivalent to projections of certain Lie algebra root systems.

The enumeration of ad-nilpotent and abelian ideals of Borel subalgebras of basic classical Lie superalgebras is given in Tables~\ref{Ta:enum_nonexcep} and~\ref{Ta:enum_excep}. This can also be achieved for parabolic subalgebras of basic classical Lie superalgebras, although we do not include this here as it can be determined when needed using the methods in this paper.

We outline the contents of this paper. In~\S\ref{s:para} we review ideals of parabolic subalgebra of Lie (super)algebras and their connection to Kostant root systems. In~\S\ref{s:conn} we prove Theorem~\ref{T:bij} by linking the super and non-super theories via Dimitrov and Fioresi's work. In~\S\ref{S:cover} we define \textit{standard covers} which describe parabolic subalgebras of $\g$ that cover all Borel subalgebras of $\G$ with respect to Theorem~\ref{T:bij}. We find standard covers of basic classical Lie superalgebras in Theorem~\ref{T:cover} by investigating associated parabolic subalgebras. In~\S\ref{s:enumerate} we apply Righi's work to these covers which leads to Tables~\ref{Ta:enum_nonexcep} and~\ref{Ta:enum_excep}.

\section{Ad-Nilpotent Ideals of Parabolic Subalgebras}\label{s:para}

\subsection{Parabolic Subalgebras of Lie Algebras}\label{ss:lie_para}
Let $\g$ be a finite-dimensional semisimple Lie algebra over $\C$. We refer the reader to~\cite{BG,DF1,DR,Ko2} as general references on the proceeding structures. Let $\p$ be a parabolic subalgebra of $\g$, so $\p$ contains a Borel subalgebra $\b=\h\oplus\n$ of $\g$ for $\h$ a Cartan subalgebra. Let $\g_{\alpha}$ for $\alpha\in\h^*$ be the root spaces of $\g$, $\Delta\subseteq\h^*$ the root system of $\g$, $\Pi$ the base of $\Delta$ that corresponds to $\b$, and $\Delta(\Pi)\coloneqq(\Z_{\geq 0}\Pi)\cap\Delta$ the resulting set of positive roots. The Levi decomposition of $\p$ is then $\p=\m\oplus\r$, where
\begin{align*}
\m=\h\oplus\bigoplus_{\alpha\in\Delta_S}\g_{\alpha},&&\r=\bigoplus_{\alpha\in\Delta(\Pi)\backslash\Delta_S}\g_{\alpha},
\end{align*}
for $S\subseteq\Pi$ and $\Delta_S\coloneqq(\Z S)\cap\Delta$ the subsystem of $\Delta$ spanned by $S$.

The center $\t=\z(\m)$ of the Levi factor $\m$ of $\p$ determines the $\t$-root spaces $\g_{\mu}\coloneqq\{x\in\g:[t,x]=\mu(t)x\text{ for all }t\in\t\}$ for $\mu\in\t^*$, $\t$-root system $\Delta^S\coloneqq\{\mu\in\t^*\backslash\{0\}:\g_{\mu}\neq 0\}$, and $\t$-root space decomposition $\g=\m\oplus\bigoplus_{\mu\in \Delta^S}\g_{\mu}$. The restriction $\pi:\h^*\rightarrow\t^*$ induces the descriptions
\begin{align*}
\Delta^S\coloneqq\pi(\Delta)\backslash\{0\},&&\Delta_S=\Delta\cap\ker(\pi),&&\g_{\mu}=\bigoplus_{\alpha\in\pi^{-1}(\mu)}\g_{\alpha},
\end{align*}
for $\mu\in\Delta^S$. We call $\Delta^S$ a \textit{Kostant root system} as named by Dimitrov and Roth in~\cite{DR}. The $\Delta^S$ notation is inspired by Brundan and Goodwin in~\cite{BG}.

Sets of positive roots of $\Delta^S$ are defined as in~\cite{DR} to be subsets $(\Delta^S)^+\subseteq\Delta^S$ where exactly one of $\mu$ and $-\mu$ is in $(\Delta^S)^+$ for each $\mu\in\Delta^S$, and where $\mu,\nu\in(\Delta^S)^+$ and $\mu+\nu\in\Delta^S$ implies $\mu+\nu\in(\Delta^S)^+$. It follows that $\Delta^S(\Pi)\coloneqq\pi(\Delta(\Pi))\backslash\{0\}$ is a set of positive roots of $\Delta^S$, and that $\r$ decomposes as the direct sum of $\t$-root spaces $\r=\bigoplus_{\mu\in\Delta^S(\Pi)}\g_{\mu}$.

As $\p$ and $\Pi$ identify the same set of positive roots in $\Delta^S$, we define $\Delta^S_{\p}\coloneqq\Delta^S(\Pi)$.

An ad-nilpotent ideal $\i$ of $\p$ is an ideal of $\p$ where each $x\in\i$ is ad-nilpotent in $\p$. These are equivalently the ideals of $\p$ contained in $\r$. A key result proved by Kostant is~\cite[Theorem 0.1]{Ko2}, which states that $\g_{\mu}$ is an irreducible \admod{$\m$} for $\mu\in\Delta^S$, and that every irreducible \adsubmod{$\m$} of $\r$ is of this form. Furthermore, $[\g_{\mu},\g_{\nu}]=\g_{\mu+\nu}$ for $\mu,\nu\in \Delta^S$ with $\mu\neq-\nu$ by~\cite[Theorem 2.3]{Ko2}. As $\h$ acts diagonally on $\g$, it follows that there exists a subset of positive roots $\Delta^S_{\i}\subseteq\Delta^S_{\p}$ such that
\begin{equation}\label{e:ideal}
\i=\bigoplus_{\mu\in\Delta^S_{\i}}\g_{\mu}.
\end{equation}
This leads to a chain of implications:
\begin{align}
[\p,\i]\subseteq \i&&\text{implies}&&\sum_{\substack{\mu\in\Delta^S_{\p} \\ \nu\in\Delta^S_{\i}}}[\g_{\mu},\g_{\nu}]\subseteq \i&&\text{implies}&&\sum_{\substack{\mu\in\Delta^S_{\p} \\ \nu\in\Delta^S_{\i}}}\g_{\mu+\nu}\subseteq\bigoplus_{\rho\in\Delta^S_{\i}}\g_{\rho}.\label{e:root_ideal}
\end{align}
There are some additional implications if $\i$ is abelian:
\begin{align}
[\i,\i]=0
&&\text{implies}&&\sum_{\mu,\nu\in\Delta^S_{\i}}[\g_{\mu},\g_{\nu}]=0&&\text{implies}&&\sum_{\mu,\nu\in\Delta^S_{\i}}\g_{\mu+\nu}=0.\label{e:root_abelian}
\end{align}
The results of these implications motivate the following definitions for subsets $I\subseteq\Delta^S_{\p}$.

\begin{enumerate}

\item[(1)] We call $I$ an \textit{ideal} of $\Delta^S_{\p}$ if $\alpha\in\Delta^S_{\p}$, $\beta\in I$, and $\alpha+\beta\in\Delta^S_{\p}$ implies $\alpha+\beta\in I$. The set of ideals of $\Delta^S_{\p}$ is denoted $\Ad(\Delta^S_{\p})$.

\item[(2)] We call $I$ an \textit{abelian ideal} of $\Delta^S_{\p}$ if $I\in \Ad(\Delta^S_{\p})$, and $\alpha,\beta\in I$ implies $\alpha+\beta\notin\Delta^S_{\p}$. The set of abelian ideals of $\Delta^S_{\p}$ is denoted $\Ab(\Delta^S_{\p})$.
\end{enumerate}

Let $\Ad(\p)$ and $\Ab(\p)$ denote the set of ad-nilpotent and abelian ideals of $\p$ respectively. The above arguments show that $\i\mapsto\Delta^S_{\i}$ is a bijection from $\Ad(\p)$ to $\Ad(\Delta^S_{\p})$, and $\Ab(\p)$ to $\Ab(\Delta^S_{\p})$. The inverse in each case takes the form $I\mapsto\bigoplus_{\mu\in I}\g_{\mu}$. It is this root correspondence that Peterson and Righi exploited in their respective proofs.

\subsection{Parabolic Subalgebras of Lie Superalgebras}\label{ss:super_para}
We refer the reader to~\cite{CW,Mu} as general references on Lie superalgebras. Of the simple Lie superalgebras classified in the 1970s by Kac in~\cite{Ka}, basic classical Lie superalgebras has the most in common with simple Lie algebras. The basic classical Lie superalgebras that are not Lie algebras are as follows.
\begin{equation*}
\overbrace{
\underbrace{\begin{matrix}
\sl(m|n) & \psl(n|n) & \osp(M|2n)
\\
m\neq n\geq 1 & n\geq 2 & M,n\geq 1
\end{matrix}}_{\text{non-exceptional}}
\hspace{0.3cm}
\underbrace{\begin{matrix}
\rD(2,1;a) & \rF(4) & \rG(3).
\\
a\in\C\backslash\{-1,0\} & &
\end{matrix}}_{\text{exceptional}}
}^{\text{basic classical}}
\end{equation*}
Our work relies only on the root system as a set and the addition relations between roots. The choice of $a\in\C\backslash\{-1,0\}$ for $\rD(2,1;a)$ thus has no affect in this paper.

The non-exceptional basic classical Lie superalgebras are reclassified based on properties of their root systems. This is presented in Table~\ref{Tb:abcd} for ease of reference. 

\begin{table}[ht]
\centering
\begin{tblr}{|c|c|c|}
\hline
Type & Lie Superalgebra & Conditions \\
\hline 
$\rA(m,n)$ & $\sl(m+1|n+1)$ & $m>n\geq 0$\\\hline 
$\rA(n,n)$ & $\psl(n+1|n+1)$ & $n\geq 1$\\\hline 
$\rB(m,n)$ & $\osp(2m+1|2n)$ & $m\geq 0$, $n\geq 1$\\\hline 
$\rC(n)$ & $\osp(2|2n-2)$ & $n\geq 2$\\\hline 
$\rD(m,n)$ & $\osp(2m|2n)$ & $m\geq 2$, $n\geq 1$\\
\hline 
\end{tblr}
\caption{The non-exceptional basic classical Lie superalgebras.}
\label{Tb:abcd}
\end{table}

Let $\G=\G\0\oplus\G\1$ be a basic classical Lie superalgebra excluding $\psl(n|n)$, $\h$ a Cartan subalgebra of $\G\0$, and $\Phi\subseteq\h^*$ the root system of $\G$ with root spaces $\G_{\alpha}$. Borel subalgebras of Lie superalgebras are not maximally solvable subalgebras in general, so alternate definitions are used. In this paper, a Borel subalgebra $\B$ of $\G$ is a subalgebra of the form $\B=\h\oplus\n$, where $\n=\bigoplus_{\alpha\in\Phi^+}\G_{\alpha}$ for a set of positive roots $\Phi^+$ of $\Phi$ defined as in~\S\ref{ss:lie_para}. There exists a base $\Sigma\subseteq\Phi$ that corresponds to $\Phi^+$, where $\R\Sigma=\Phi$, and $\alpha\in\Phi^+$ uniquely decomposes as $\alpha=\sum_{\beta\in\Sigma}k_{\beta}\beta$ with $k_{\beta}\in\Z_{\geq 0}$ for $\beta\in\Sigma$.

Let $\P$ a parabolic subalgebra of $\G$, so $\P$ contains a Borel subalgebra $\B$ of $\G$. There exists a decomposition $\P=\m\oplus\r$ with $\m$, $\r$, and $\Phi_S$ for $S\subseteq\Sigma$ defined as in~\S\ref{ss:lie_para}. We refer the reader to~\cite{DF1} as a general reference on Kostant root systems of Lie superalgebras. Many of the concepts outlined~\S\ref{ss:lie_para} extend appropriately to Lie superalgebras. The center $\t=\z(\m)$ determines the $\t$-root spaces $\G_{\mu}$, $\t$-root system or Kostant root system $\Phi^S$, and $\t$-root space decomposition $\G=\m\oplus\bigoplus_{\mu\in\Phi^S}\G_{\mu}$. The restriction $\pi:\h^*\rightarrow\t^*$ gives comparable forms for $\Phi_S$, $\Phi^S$, and $\G_{\mu}$. We define $\Phi^S_{\P}$ similarly.

An ideal $\i$ of $\P$ is a $\Z_2$-graded subspace of $\P$ where $[\P,\i]\subseteq\i$. Ad-nilpotent ideals of $\P$ are defined like for $\p$, and are the ideals of $\P$ contained in $\r$ for the same reasons. As stated in~\cite[Theorem 1.4]{DF1}, each $\t$-root space $\G_{\mu}$ is an irreducible \admod{$\m$}, and $[\G_{\mu},\G_{\nu}]=\G_{\mu+\nu}$ for $\mu,\nu\in\Phi^S$ with $\mu\neq-\nu$. It follows that an ad-nilpotent ideal $\i$ of $\P$ decomposes as in~\eqref{e:ideal}, so there exists a subset of positive roots $\Phi^S_{\i}\subseteq\Phi^S_{\P}$ such that
\begin{equation}\label{e:ideal_s}
\i=\bigoplus_{\mu\in\Phi^S_{\i}}\G_{\mu}.
\end{equation}
Implications like those in~\eqref{e:root_ideal} and~\eqref{e:root_abelian} also hold.

We define ideals and abelian ideals of $\Phi^S_{\P}$, and $\Ad(\Phi^S_{\P})$, $\Ab(\Phi^S_{\P})$, $\Ad(\P)$, and $\Ab(\P)$ as in~\S\ref{ss:lie_para}. We again have the bijection $\i\mapsto\Phi^S_{\i}$ from $\Ad(\P)$ to $\Ad(\Phi^S_{\P})$, and $\Ab(\P)$ to $\Ab(\Phi^S_{\P})$.

\begin{remark}\label{R:pslnn}
The above decomposition exists for ad-nilpotent and abelian ideals of Borel subalgebras $\B$ of $\psl(n|n)$ for $n>2$ due to~\cite[Theorem 1.18(2)]{CW} and~\cite[Lemma 2.1.1]{Mu}. This implies the same bijection from $\Ad(\B)$ to $\Ad(\Phi_{\B})$, and $\Ab(\B)$ to $\Ab(\Phi_{\B})$. This decomposition fails for $\psl(2|2)$, which is the only basic classical Lie superalgebras to have 2-dimensional root spaces. Consequently, Borel subalgebras of $\psl(2|2)$ have infinitely many ad-nilpotent ideals. It is for these reasons that we do not consider $\psl(n|n)$ again until~\S\ref{s:enumerate}.
\end{remark}

\section{The Connection}\label{s:conn}
In this section we prove Theorem~\ref{T:bij} by relating $\p$ and $\P$ using results in~\cite{DF1}. The cornerstone of our approach is the proof of Theorem 1.4 of \textit{loc. cit.}, where it was shown that for a root system $\Phi$ of a basic classical Lie superalgebra $\G$ excluding $\psl(n|n)$, there exists a root system $\Delta$ of a simple Lie algebra $\g$, subset $S\subsetneq\Pi$ of a base, and vector space isomorphism $Q:\R\Delta^S\rightarrow\R\Phi$ where $Q(\Delta^S)=\Phi$. We call $\Delta^S$ and $\Phi$ \textit{equivalent} in this case. 

There are multiple such subsets used throughout~\cite{DF1}, with many appearing in~\cite{DR} and~\cite{DF2}. Table~\ref{Ta:quotients} lists the subsets described in~\cite[Theorem 4.7]{DF2}, except when $\G=\rD(2,1;a)$, where an equivalent choice is made to simplify future calculations. The subsets given in~\textit{loc. cit.} are of the form $\Pi\backslash S$ for $\Pi$ the base given in~\cite[\S 4]{Bo}.

\begin{table}[ht]
\centering
\begin{tblr}{|c|c|c|c|}
\hline[1.5pt]
$\G$ & $\g$ & $\Delta$ & $S$
\\
\hline[1.5pt]
$\rA(m,n)$, $m\neq n$ & $\sl_{m+n+2}(\C)$ & $A_{m+n+1}$ & $\O$
\\
\hline
$\rB(m,n)$ & $\so_{2(m+2n)+1}(\C)$ & $B_{m+2n}$ & $\{\alpha_{m+1},\alpha_{m+3},\ldots,\alpha_{m+2n-1}\}$
\\
\hline
$\rC(n)$ & $\so_{2(2n-1)}(\C)$ & $D_{2n-1}$ & $\{\alpha_{2},\alpha_{4},\ldots,\alpha_{2n-2}\}$
\\
\hline
$\rD(m,n)$ & $\so_{2(m+2n)}(\C)$ & $D_{m+2n}$ & $\{\alpha_{m+1},\alpha_{m+3},\ldots,\alpha_{m+2n-1}\}$
\\
\hline
$\rD(2,1;a)$ & $\so_{8}(\C)$ & $D_{4}$ & $\{\alpha_1\}$
\\
\hline
$\rF(4)$ & $\e_{7}$ & $E_7$ & $\{\alpha_2,\alpha_3,\alpha_4\}$
\\
\hline
$\rG(3)$ & $\e_{8}$ & $E_8$ & $\{\alpha_2,\alpha_4,\alpha_5,\alpha_6,\alpha_7\}$
\\
\hline
\end{tblr}
\caption{(left to right) Lie superalgebra $\G$; Lie algebra $\g$; root system $\Delta$ of $\g$; subset $S\subsetneq\Pi$ such that $\Delta^S$ is equivalent to $\Phi$.}
\label{Ta:quotients}
\end{table}

We give an example of how this can be shown.

\begin{example}\label{Ex:equiv_system}
Let $\G=\rD(2,1;a)$. The root system of $\G$ as denoted by row vectors in~\cite[\S 4.2]{Mu} is $\Phi=\{(\pm 2,0,0),(0,\pm 2,0),(0,0,\pm 2),(\pm1,\pm1,\pm1)\}$.

Let $\g=\so_{8}(\C)$. The root system of $\g$ as denoted in~\cite[\S 4]{Bo} is $\Delta=\{\epsilon_i\pm \epsilon_j:1\leq i\neq j\leq 4\}$. The base $\Pi=\{\alpha_1,\alpha_2,\alpha_3,\alpha_4\}$ from \textit{loc. cit.} has elements of the form
\begin{align*}
\alpha_1=\epsilon_1-\epsilon_2,&&\alpha_2=\epsilon_2-\epsilon_3,&&\alpha_3=\epsilon_3-\epsilon_4,&&\alpha_4=\epsilon_3+\epsilon_4.
\end{align*}
To show that $\Phi$ and $\Delta^{\{\alpha_1\}}$ are equivalent as stated in Table~\ref{Ta:quotients}, it is sufficient to show that there exists a vector space isomorphism $Q:\R\Delta^{\{\alpha_1\}}\rightarrow\R\Phi$ such that $Q(\Delta^{\{\alpha_1\}}(\Pi))$ is a set of positive roots of $\Phi$. We construct $\Delta^{\{\alpha_1\}}(\Pi)$ in the below figure by grouping the preimage $\Delta(\Pi)\cap\pi^{-1}(\mu)$ of each $\mu\in\Delta^{\{\alpha_1\}}(\Pi)$. The shaded square indicates that $\pi(\epsilon_1-\epsilon_2)=0$.

\bigskip
\begin{center}
\includegraphics[scale=0.5]{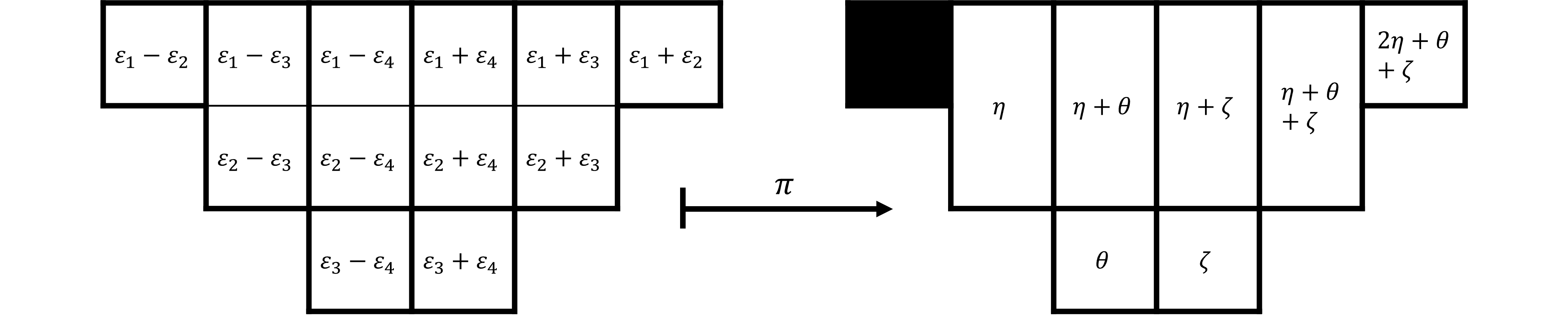}
\end{center}

\noindent
Consider both the set of positive roots
\begin{equation*}
\Phi^+=\{(1,-1,-1),(0,2,0),(0,0,2),(1,1,-1),(1,-1,1),(1,1,1),(2,0,0)\},
\end{equation*}
and the map $Q:\R\Delta^{\{\alpha_1\}}\rightarrow\R\Phi$, given by
\begin{align*}
Q(\eta)&=(1,-1,-1),&Q(\theta)&=(0,2,0),
&Q(\zeta)&=(0,0,2),&
\end{align*}
and extended linearly. One may verify that $Q$ is a bijection with $Q(\Delta^{\{\alpha_1\}}(\Pi))=\Phi^+$, and hence that $\Phi$ and $\Delta^{\{\alpha_1\}}$ are equivalent. This similarly works for $\{\alpha_2\}$, $\{\alpha_3\}$, and $\{\alpha_4\}$.
\end{example}

The following definition is motivated by the method used in Example~\ref{Ex:equiv_system}.

\begin{definition}\label{D:analogous}
Parabolic subalgebras $\p$ and $\P$, with sets of positive roots $\Delta^R_{\p}$ and $\Phi^S_{\P}$, are \textit{analogous} if there exists a vector space isomorphism $Q:\R\Delta^R_{\p}\rightarrow\R\Phi^S_{\P}$ where $Q(\Delta^R_{\p})=\Phi^S_{\P}$.
\end{definition}

Preserving addition relations between roots, the vector space isomorphism $Q$ naturally induces a bijection from $\Ad(\Delta^R_{\p})$ to $\Ad(\Phi^S_{\P})$, and $\Ab(\Delta^R_{\p})$ to $\Ab(\Phi^S_{\P})$. The following lemma is an immediate consequence of this and the ideal bijections outlined in~\S\ref{s:para}.

\begin{lemma}\label{L:para}
Let $\p$ and $\P$ be analogous parabolic subalgebras of $\g$ and $\G$ respectively with $\G\neq\psl(n|n)$. Then both $\Ad(\p)$ and $\Ad(\P)$, and $\Ab(\p)$ and $\Ab(\P)$, are in bijection.
\end{lemma}

We are now in a position to prove our main theorem.

\begin{proof}[Proof of Theorem~\ref{T:bij}]
Fix $\g$ from the same row of Table~\ref{Ta:quotients} as $\G$. As discussed at the beginning of~\S\ref{s:conn}, there exists $S\subseteq\Pi$ and $Q:\R\Delta^S\rightarrow\R\Phi$ such that $Q(\Delta^S)=\Phi$.

We claim that there exists a parabolic subalgebra $\p$ of $\g$ that is analogous to $\P$.

Observe from~\cite[Proposition 3.4]{DF1} that
\begin{align}\label{e:positive_bijection}
\{\text{bases of }\Delta\text{ that contain }S\}
\rightarrow
\{\text{bases of }\Delta^S\}
\rightarrow
\{\text{bases of }\Phi\}
\end{align}
is a chain of bijections given by $B\mapsto\pi(B)\backslash\{0\}\mapsto Q(\pi(B)\backslash\{0\})$. Therefore there exists a base $B\supseteq S$ of $\Delta$, and subset $U\subsetneq Q(\pi(B)\backslash\{0\})$, such that $\Phi^U$ is the Kostant root system determined by $\P$. It follows that $\Delta^R$, for $R=\Delta\cap(Q\circ\pi)^{-1}(U)$, is equivalent to $\Phi^U$.

Consider a parabolic subalgebra $\p=\m\oplus\r$ of $\g$ with $\r=\bigoplus_{\mu\in\Delta^R(B)}\g_{\mu}$. This is analogous to $\P$, proving our claim. The result follows from Lemma~\ref{L:para}.
\end{proof}

\section{Standard Covers}\label{S:cover}

Fix $\G$, $\g$, and $S$ from a row of Table~\ref{Ta:quotients} throughout. Given Theorem~\ref{T:bij}, the enumeration of ad-nilpotent (respectively, abelian) ideals of parabolic subalgebras of $\G$ is equal to that of certain parabolic subalgebras of $\g$. We now determine these parabolic subalgebras of $\g$.

We may consider only standard parabolic subalgebras of $\g$, as each parabolic subalgebra of $\g$ is isomorphic to a standard parabolic subalgebra. Let the base $\Pi$ of $\Delta$ that corresponds to the fixed Borel subalgebra that they contain be as in~\cite[\S 4]{Bo}, and denote by $\p_S$ the standard parabolic subalgebra that determines the Kostant root system $\Delta^S$ for $S\subseteq\Pi$.

Recall that $\p_R$ and $\p_S$ are \textit{associated} if their Levi factors are conjugate under the Weyl group $W$ of $\Delta$, or equivalently if $w(R)=S$ for some $w\in W$.

These concepts motivate the below lemma.

\begin{lemma}\label{L:cover_assoc}
If $\p_R$ is associated to $\p_S$, then $\p_R$ is analogous to a Borel subalgebra of $\G$. And for all Borel subalgebras $\B$ of $\G$, there exists $R\subsetneq\Pi$ such that $\p_R$ is associated to $\p_{S}$ and analogous to $\B$.
\end{lemma}

\begin{proof}
The result immediately holds if $S=\O$, so suppose $S\neq\O$.

Let $\p_R$ be associated to $\p_S$, so $w(R)=S$ for some $w\in W$. Note that $w(\Pi)\supsetneq S$ is a base of $\Delta$. A parabolic subalgebra $\p=\m\oplus\r$ of $\g$ with $\r=\bigoplus_{\mu\in\Delta^{S}(w(\Pi))}\g_{\mu}$ is thus analogous to some Borel subalgebra $\B$ of $\G$ by~\eqref{e:positive_bijection}. Hence $\p_R$, being conjugate to $\p$, is analogous to $\B$.

Let $\B$ be a Borel subalgebra of $\G$. From the proof of Theorem~\ref{T:bij} there exists a base $B\supseteq S$ of $\Delta$, where a parabolic subalgebra $\p=\m\oplus\r$ of $\g$ with $\r=\bigoplus_{\mu\in\Delta^{S}(B))}\g_{\mu}$ is analogous $\B$. Let $w\in W$ be such that $w(B)=\Pi$. It follows that $\p_{w(S)}$ is analogous to $\B$.
\end{proof}

Let $\P$ be a parabolic subalgebra of $\G$ containing a Borel subalgebra $\B$. It first follows from Lemma~\ref{L:cover_assoc} that $\p_R$ is analogous to $\B$ for some $R\subsetneq\Pi$, and then from the proof of Theorem~\ref{T:bij} that $\p_{R'}$ is analogous to $\P$ for some $S\subseteq R'\subsetneq\Pi$. Standard Parabolic subalgebras of $\g$ associated to $\p_S$ are thus sufficient for studying parabolic subalgebras of $\G$. Note also that if $\p_R$ is analogous to $\B$, then $\p_{g(R)}$ is analogous to $\B$ for $g$ an automorphism of the Dynkin diagram of $\Delta$ with respect to $\Pi$. We define the following with this in mind.

\begin{definition}
A \textit{standard cover} of $\G$ over $\g$ is a subset $\cP\subseteq\{w(S)\subseteq\Pi:w\in W\}$, which is maximal in the sense that $g(R)\neq R'$ for distinct $R,R'\in\cP$ and automorphisms $g$ of the Dynkin diagram of $\Delta$ with respect to $\Pi$.
\end{definition}

The $\h$-stable subalgebras of $\g$ were classified in the 1950s by Dynkin in~\cite[p.146-148]{Dy}. Dynkin's work shows, with a few exceptions, that $\p_R$ and $\p_S$ are associated if and only if $R$ and $S$ have isomorphic Dynkin subdiagrams where root lengths are taken into account. The only exception that applies to $S$ of Table~\ref{Ta:quotients} is when $\Delta=D_l$. In this case, $\p_R$ and $\p_S$ are associated if and only if $R$ and $S$ have isomorphic Dynkin subdiagrams and $R\not\supseteq\{\alpha_{l},\alpha_{l-1}\}$.

We use the above results to find standard covers of basic classical Lie superalgebras. The following Dynkin diagrams are labelled with respect to $\Pi$.

\begin{theorem}\label{T:cover}
There exists a standard cover $\cP$ of $\G$ over $\g$ as follows:
\begin{enumerate}

\item[(i)] if $\G=\rA(m,n)$ and $\g=\sl_{m+n+2}(\C)$ for $m\neq n$, then
\begin{equation*}
\cP=\{\O\};
\end{equation*}

\item[(ii)] if $\G=\rB(m,n)$ and $\g=\so_{2(m+2n)+1}(\C)$, then 
\begin{equation*}
\cP=\{\{\alpha_{i_j}\}_{j=1}^n\subsetneq\Pi\backslash\{\alpha_{m+2n}\}:|i_{j}-i_{k}|>1\text{ for all }1\leq j< k\leq n-1\};
\end{equation*}

\item[(iii)] if $\G=\rC(n)$ and $\g=\so_{2(2n-1)}(\C)$, then
\begin{equation*}
\cP=\{\{\alpha_{2j-1}\}_{j=1}^c\cup\{\alpha_{2j}\}_{j=c+1}^{n-1}\subsetneq\Pi:c=0,\ldots,n-1\};
\end{equation*}

\item[(iv)] if $\G=\rD(m,n)$ and $\g=\so_{2(m+2n)}(\C)$, then
\begin{equation*}
\cP=\{\{\alpha_{i_j}\}_{j=1}^n\subsetneq\Pi\backslash\{\alpha_{m+2n}\}:|i_{j}-i_{k}|>1\text{ for all }1\leq j< k\leq n-1\}
\end{equation*}
for $(m,n)\neq(2,1)$, and $\cP=\{\{\alpha_1\},\{\alpha_2\}\}$ for $(m,n)=(2,1)$;

\item[(v)] if $\G=\rD(2,1;a)$ and $\g=\so_{8}(\C)$, then
\begin{equation*}
\cP=\{\{\alpha_1\},\{\alpha_2\}\};
\end{equation*}

\item[(vi)] if $\G=\rF(4)$ and $\g=\e_{7}$, then
\begin{equation*}
\cP=\left\{\begin{matrix}
\{\alpha_{1},\alpha_{3},\alpha_{4}\},\{\alpha_{2},\alpha_{3},\alpha_{4}\},\{\alpha_{2},\alpha_{4},\alpha_{5}\},
\\
\{\alpha_{3},\alpha_{4},\alpha_{5}\},\{\alpha_{4},\alpha_{5},\alpha_{6}\},\{\alpha_{5},\alpha_{6},\alpha_{7}\}
\end{matrix}\right\};
\end{equation*}

\item[(vii)] if $\G=\rG(3)$ and $\g=\e_{8}$, then
\begin{equation*}
\cP=\left\{\begin{matrix}
\{\alpha_{1},\alpha_{3},\alpha_{4},\alpha_{5},\alpha_{6}\},\{\alpha_{2},\alpha_{4},\alpha_{5},\alpha_{6},\alpha_{7}\},
\\
\{\alpha_{3},\alpha_{4},\alpha_{5},\alpha_{6},\alpha_{7}\}$, $\{\alpha_{4},\alpha_{5},\alpha_{6},\alpha_{7},\alpha_{8}\}\end{matrix}\right\}.
\end{equation*}
\end{enumerate}
\end{theorem}

\begin{proof}
(i) This immediately holds as $S=\O$ in this case.

\noindent
(ii) View the Dynkin subdiagram of $S=\{\alpha_{m+1},\alpha_{m+3},\ldots,\alpha_{l-1}\}$ in $B_{l}$ for $l=m+2n$.
\begin{align*}
\begin{dynkinDiagram}[labels={1,2,m,m+1,m+2,m+3,m+4,l-2,l-1,l}, edge length=1cm, root
radius=0.1cm]{B}{oo.o*o*o.o*o}
\end{dynkinDiagram}
\end{align*}
As $\alpha_l\notin S$, we find that $\p_R$ and $\p_S$ are associated if and only if $\# R=n$, $\alpha_{l}\notin R$, and all roots in $R$ are mutually non-adjacent. The result follows.

\noindent
(iii),(iv) View the Dynkin subdiagram of $S=\{\alpha_{m+1},\alpha_{m+3},\ldots,\alpha_{l-1}\}$ in $D_{l}$ for $l=m+2n$.
\begin{align*}
\begin{dynkinDiagram}[labels={1,2,m,m+1,m+2,m+3,m+4,l-4,l-3,,l-1,l}, labels*={,,,,,,,,,l-2,}, edge length=1cm, root
radius=0.1cm]{D}{oo.o*o*o.o*o*o}
\end{dynkinDiagram}
\end{align*}
We find that $\p_R$ and $\p_S$ are associated if and only if $\# R=n$, $R\not\supseteq\{\alpha_{l},\alpha_{l-1}\}$, and all roots in $R$ are mutually non-adjacent. Due to the symmetries of $D_l$, we take $\cP$ to contain all such subsets with $\alpha_{l}\notin R$. The result follows for $\rD(m,n)$ with $(m,n)\neq(2,1)$, and for $\rC(n+1)$ by setting $m=1$. The $\rD(2,1)$ case follows from (v) of this theorem.

\noindent
(v) View the Dynkin subdiagram of $S=\{\alpha_{1}\}$ in $D_{4}$.
\begin{align*}
\begin{dynkinDiagram}[labels={1,2,3,4}, edge length=1cm, root
radius=0.1cm]{D}{*ooo}
\end{dynkinDiagram}
\end{align*}
We find that $\p_R$ and $\p_S$ are associated if and only if $\# R=1$. Due to the symmetries of $D_4$, we take $\cP=\{\{\alpha_1\},\{\alpha_2\}\}$.

\noindent
(vi) View the Dynkin subdiagram of $S=\{\alpha_2,\alpha_3,\alpha_4\}$ in $E_{7}$.
\begin{align*}
\begin{dynkinDiagram}[labels={1,2,3,4,5,6,7}, edge length=1cm, root
radius=0.1cm]{E}{o***ooo}
\end{dynkinDiagram}
\end{align*}
We find that $\p_R$ and $\p_S$ are associated if and only if $R$ has an $A_3$ type Dynkin subdiagram. These correspond to the below sets of shaded roots.
\begin{align*}
\begin{dynkinDiagram}[edge length=0.5cm, root
radius=0.1cm]{E}{*o**ooo}
\end{dynkinDiagram},&&\begin{dynkinDiagram}[edge length=0.5cm, root
radius=0.1cm]{E}{o***ooo}
\end{dynkinDiagram},&&\begin{dynkinDiagram}[edge length=0.5cm, root
radius=0.1cm]{E}{o*o**oo}
\end{dynkinDiagram},
\\
\begin{dynkinDiagram}[edge length=0.5cm, root
radius=0.1cm]{E}{oo***oo}
\end{dynkinDiagram},&&\begin{dynkinDiagram}[edge length=0.5cm, root
radius=0.1cm]{E}{ooo***o}
\end{dynkinDiagram},&&\begin{dynkinDiagram}[edge length=0.5cm, root
radius=0.1cm]{E}{oooo***}
\end{dynkinDiagram}.
\end{align*}

\noindent
(vii) View the Dynkin subdiagram of $S=\{\alpha_2,\alpha_4,\alpha_5,\alpha_6,\alpha_7\}$ in $E_{8}$.
\begin{align*}
\begin{dynkinDiagram}[labels={1,2,3,4,5,6,7,8}, edge length=1cm, root
radius=0.1cm]{E}{o*o****o}
\end{dynkinDiagram}
\end{align*}
We find that $\p_R$ and $\p_S$ are associated if and only if $R$ has an $A_5$ type Dynkin subdiagram. These correspond to the below sets of shaded roots.
\begin{align*}
\begin{dynkinDiagram}[edge length=0.5cm, root
radius=0.1cm]{E}{*o****oo}
\end{dynkinDiagram},&&\begin{dynkinDiagram}[edge length=0.5cm, root
radius=0.1cm]{E}{o*o****o}
\end{dynkinDiagram},
\\
\begin{dynkinDiagram}[edge length=0.5cm, root
radius=0.1cm]{E}{oo*****o}
\end{dynkinDiagram},&&\begin{dynkinDiagram}[edge length=0.5cm, root
radius=0.1cm]{E}{ooo*****}
\end{dynkinDiagram}.
\end{align*}
This concludes the proof.
\end{proof}

This author speculates that the standard covers $\cP$ of Theorem~\ref{T:cover}, except the $\G=\rC(2)$ case, are minimal in the sense that for all Borel subalgebras $\B$ of $\G$, there exists a unique $S\in\cP$ where $\p_S$ is analogous $\B$. This holds for $\rA(m,n)$ as $\cP=\{\O\}$, and for $\rD(2,1;a)$, $\rF(4)$, and $\rG(3)$ as Table~\ref{Ta:enum_excep} contains no duplicate columns. It remains to be shown for the orthosymplectic cases. This question of minimality does not affect the results of~\S\ref{s:enumerate}.

For each Borel subalgebra $\B$ of $\G$, an element $S\in\cP$ where $\p_S$ is analogous to $\B$ may be resolved in a standard way. For example if $\G=\rB(m,n)$ and $\g=\so_{2(m+2n)+1}(\C)$, then the non-isotropic odd roots $\alpha\in\Delta_{\B}$ where $2\alpha\in\Delta_{\B}$ completely determine $S$. This is because $\alpha_i=\epsilon_i-\epsilon_{i+1}\in S$ induces the Kostant roots $\mu,2\mu\in\Delta^S$ with root spaces $\g_{\mu}=\g_{\epsilon_i}\oplus\g_{\epsilon_{i+1}}$ and $\g_{2\mu}=\g_{\epsilon_i+\epsilon_{i+1}}$. We do not consider the procedures in each case for technical reasons.

\section{Enumeration of Ad-Nilpotent Ideals of Borel Subalgebras}\label{s:enumerate}
We are now in a position to enumerate ad-nilpotent and abelian ideals of Borel subalgebras of basic classical Lie superalgebras. For $\G\neq\psl(n|n)$, this is done by substituting $S\in\cP$ from Theorem~\ref{T:cover} into the appropriate result in~\cite{Ri1,
Ri2}. We then implicitly use Lemma~\ref{L:para}.

We do not consider $\rA(1,1)$, $\rC(2)$, or $\rD(2,1)$ as the $\rA(1,1)$ case is explained in Remark~\ref{R:pslnn}, $\rC(2)\cong\rA(1,0)$, and $\rD(2,1)$ and $\rD(2,1;a)$ have equivalent root systems.

Let $\G=\psl(n|n)$, $\g=\sl_{2n}(\C)$, and $n>2$. Despite $\rank(\G)=\rank(\g)-1$, the root systems of $\G$ and $\g$ are equivalent as sets and have comparable addition relations between roots. Given Remark~\ref{R:pslnn}, there is a natural bijection of ad-nilpotent (respectively, abelian) ideals of Borel subalgebras $\B$ of $\G$ and $\b$ of $\g$, so $\#\Ad(\B)=\#\Ad(\b)$ and $\#\Ab(\B)=\#\Ab(\b)$. These are calculated from the same results used in the $\rA(m,n)$ case of Theorem~\ref{T:enumerate}.

We fix the following conventions and notations for the remainder of this paper.
\begin{itemize}
\item Fix $\G$, $\g$, and $\cP$ from a case of Theorem~\ref{T:cover} with $\Delta$ the root system of $\g$.

\item For each Borel subalgebra $\B$ of $\G$, fix $S\in\cP$ such that $\p_{S}$ is analogous to $\B$. Label its roots $S=\{\alpha_{S_j}:j=1,\ldots,r\}$ where $i_1<\ldots<i_r$.

\item Let $c$ be as in Theorem~\ref{T:cover}(iii) for $\G=\rC(n)$, and $d=m+n$ for $\G=\rB(m,n),\rD(m,n)$.

\item Let $S_1,\ldots,S_s$ be the connected components of $S$ in the Dynkin subdiagram of $\Delta$ with respect to $\Pi$. Let $m_j=\min\{i:\alpha_i\in S_j\}$ and $r_j=\# S_j$ for $j=1,\ldots,s$. Assume $S_1,\ldots, S_s$ is indexed such that $m_1<\ldots<m_s$. 

\item Let $l_j=m_j-\sum_{k=1}^{j-1}r_k$ for $j=1,\ldots,s$.
\end{itemize}

We never have both $\Delta=D_l$ and $S_s=\{\alpha_l\}$ in Theorem~\ref{T:cover}. The forms for $l_j$ given in equations (11) and (12) of~\cite{Ri1} thus reduce to the one given.

We adopt some additional notations for the non-exceptional and $\rD(2,1;a)$ cases. Let
\begin{align*}
\cC_{l}=\frac{1}{l+1}\binom{2l}{l},&&\cT'_{p,q}=\frac{p-q+2}{p+q+2}\binom{p+q+2}{q},
\end{align*}
for $l\in\Z_{\geq 0}$ and $p,q\in\Z_{\geq 1}$ with $p\geq q$. One may verify that the equation for $\cT'_{p,q}$ is equivalent to that in~\cite[Proposition 5]{Ri1}. We consider $\cT'_{p,0}= 1$ and $\cT'_{p,q}= 0$ for $p\geq 0>q$.

\begin{table}[htp]
\centering
\scalebox{0.9}{
\begin{tblr}{|c|c|c|}
\hline[1.5pt]
  & Conditions & $\rA(m,n)$, $(m,n)\neq(1,1)$
\\
\hline[1.5pt]
$\#\Ad(\B)$ & N.A. & $\cC_{m+n+2}$
\\
\hline
$\#\Ab(\B)$ & N.A & $2^{m+n+1}$
\\
\hline[1.5pt]
  & & $\rB(m,n)$
\\
\hline[1.5pt]
$\#\Ad(\B)$ & N.A. & $(d+1)\cC_{d}+\sum_{j=1}^n\cT_{2d-i_j+j-1,i_j-j}'$
\\
\hline
\SetCell[r=2]{c} $\#\Ab(\B)$ & $\alpha_1\in S$ & $2^{d}+2\sum_{j=1}^n\binom{d-1}{i_j-j}$ 
\\
\hline
& $\alpha_1\notin S$ & $2^{d-1}+\sum_{j=1}^n\binom{d-1}{i_j-j}$
\\
\hline[1.5pt]
  & & $\rC(n)$, $n\neq 2$
\\
\hline[1.5pt]
\SetCell[r=3]{c} $\#\Ad(\B)$ & $c=0$ & $\frac{n+1}{2}\cC_{n}+\sum_{j=1}^{n-1}\cT'_{2n-j-2,j}$
\\
\hline
& \SetCell[r=2]{c} $c\in\{1,\ldots,n-1\}$ & $(3n-2)\cC_{n-1}+\sum_{j=1}^{c}\left[\cT'_{2n-j-1,j-2}+\cT'_{2n-j-1,j-1}\right]$ 
\\
&& $+\sum_{j=c+1}^{n-1}\left[\cT'_{2n-j-2,j-1}+\cT'_{2n-j-2,j}\right]$
\\
\hline
\SetCell[r=2]{c} $\#\Ab(\B)$ & $c\in\{0,n-1\}$ & $2^{n}+2^{n-1}-2$
\\
\hline
& $c\in\{1,\ldots,n-2\}$ & $2^{n}-\binom{n-1}{c}$
\\
\hline[1.5pt]
  & & $\rD(m,n)$, $(m,n)\neq(2,1)$
\\
\hline[1.5pt]
\SetCell[r=3]{c} $\#\Ad(\B)$ & $\alpha_{m+2n-1}\in S$ & $\frac{d+1}{2}\cC_{d}+\sum_{j=1}^n\cT'_{2d-i_j+j-2,i_j-j}$
\\
\hline
& \SetCell[r=2]{c} $\alpha_{m+2n-1}\notin S$ & $(3d-2)\cC_{d-1}$
\\
&& $+\sum_{j=1}^n[\cT'_{2d-i_j+j-2,i_j-j-1}+\cT'_{2d-i_j+j-2,i_j-j}]$
\\
\hline
\SetCell[r=4]{c} $\#\Ab(\B)$ & $\alpha_{1}\in S$, $\alpha_{m+2n-1}\notin S$ & $2^{d}-2^{d-2}+\sum_{j=1}^n \left[2\binom{d-1}{i_j-j}-\binom{d-2}{i_j-j}\right]$
\\
\hline
& $\alpha_{1}\notin S$, $\alpha_{m+2n-1}\in S$ & $2^{d-1}+2^{d-2}+\sum_{j=1}^n \binom{d-1}{i_j-j}+\sum_{j=1}^{n-1}\binom{d-2}{i_j-j}$
\\
\hline
& $\alpha_{1}\in S$, $\alpha_{m+2n-1}\in S$ & $2^{d-1}+\sum_{j=1}^n \binom{d-1}{i_j-j}$
\\
\hline
& $\alpha_{1}\notin S$, $\alpha_{m+2n-1}\notin S$ & $2^{d}+2\sum_{j=1}^n \binom{d-1}{i_j-j}$
\\
\hline
\end{tblr}}
\caption{The number of ad-nilpotent and abelian ideals of Borel subalgebras of non-exceptional basic classical Lie superalgebras.}
\label{Ta:enum_nonexcep}
\end{table}

\begin{table}[htp]
\scalebox{0.9}{
\begin{tblr}{|c|c|c|c|c|c|c|}
\hline[1.5pt]
  & \SetCell[c=2]{c} $\rD(2,1;a)$ && \SetCell[c=4]{c} $\rG(3)$
\\
\hline[1.5pt]
$S$ & $\begin{dynkinDiagram}[edge length=0.2cm]{D}{*ooo}
\end{dynkinDiagram}$ 
& $\begin{dynkinDiagram}[edge length=0.2cm]{D}{o*oo}
\end{dynkinDiagram}$ 
& $\begin{dynkinDiagram}[edge length=0.2cm]{E}{*o****oo}
\end{dynkinDiagram}$ 
& $\begin{dynkinDiagram}[edge length=0.2cm]{E}{o*o****o}
\end{dynkinDiagram}$ 
& $\begin{dynkinDiagram}[edge length=0.2cm]{E}{oo*****o}
\end{dynkinDiagram}$ 
& $\begin{dynkinDiagram}[edge length=0.2cm]{E}{ooo*****}
\end{dynkinDiagram}$
\\
\hline
$\#\Ad(\B)$ & $15$ & $19$ & $36$ & $34$ & $39$ & $25$
\\
\hline
$\#\Ab(\B)$ & $7$ & $12$ & $9$ & $8$ & $11$ & $11$
\\
\hline[1.5pt]
  & \SetCell[c=6]{c} $\rF(4)$
\\
\hline[1.5pt]
$S$ & $\begin{dynkinDiagram}[edge length=0.2cm]{E}{*o**ooo}
\end{dynkinDiagram}$ 
& $\begin{dynkinDiagram}[edge length=0.2cm]{E}{o***ooo}
\end{dynkinDiagram}$ 
& $\begin{dynkinDiagram}[edge length=0.2cm]{E}{o*o**oo}
\end{dynkinDiagram}$ 
& $\begin{dynkinDiagram}[edge length=0.2cm]{E}{oo***oo}
\end{dynkinDiagram}$ 
& $\begin{dynkinDiagram}[edge length=0.2cm]{E}{ooo***o}
\end{dynkinDiagram}$ 
& $\begin{dynkinDiagram}[edge length=0.2cm]{E}{oooo***}
\end{dynkinDiagram}$
\\
\hline
$\#\Ad(\B)$ & $64$ & $83$ & $80$ & $95$ & $93$ & $54$
\\
\hline
$\#\Ab(\B)$ & $20$ & $18$ & $18$ & $22$ & $20$ & $7$
\\
\hline
\end{tblr}}
\caption{The number of ad-nilpotent and abelian ideals of Borel subalgebras of exceptional basic classical Lie superalgebras.}
\label{Ta:enum_excep}
\end{table}

Coincidently, every instance of the variable $n$ in an expression from~\cite{Ri1} that we use, is the same $n$ from the appropriate case of Theorem~\ref{T:cover}. We also emphasise that~\cite{Ri1} and~\cite{Ri2} use equivalently labelled Dynkin diagrams to those in~\cite[\S 4]{Bo}, so our results are compatible.

We proceed with the enumerations.

\begin{theorem}\label{T:enumerate}
Tables~\ref{Ta:enum_nonexcep} and~\ref{Ta:enum_excep} contain the correct enumerations.
\end{theorem}

\begin{proof}
$\G=\rA(m,n)$, $(m,n)\neq (1,1)$

\noindent
The $m=n$ case has been discussed so let $m\neq n$. Note $\Delta=A_{m+n+1}$. We get $\#\Ad(\B)$ from~\cite[Theorem 5.12]{Ri1} and $\#\Ab(\B)$ from~\cite[Theorem 4.8]{Ri1}.

\vspace{0.25cm}
\noindent
$\G=\rB(m,n)$

\noindent
Note $\Delta=B_{m+2n}$, $r=s=n$, $S_j=\{\alpha_{i_j}\}$, $r_j=1$, $m_j=i_j$, and $l_j=i_j-j+1$ for $j=1,\ldots,n$. We get $\#\Ad(\B)$ from~\cite[Theorem 5.12]{Ri1} and $\#\Ab(\B)$ from~\cite[Proposition 5.18]{Ri1}.

\vspace{0.25cm}
\noindent
$\G=\rC(n)$, $n>2$

\noindent
Note $\Delta=D_{2n-1}$, $r=s=n-1$, $S_j=\{\alpha_{i_j}\}$, $r_j=1$, $m_j=i_j$, and $l_j=i_j-j+1$, where
\begin{equation*}
i_j=\begin{cases}
2j-1 &\text{if }j=1,\ldots,c,
\\
2j &\text{if }j=c+1,\ldots,n-1,
\end{cases}
\end{equation*}
for $j=1,\ldots,n-1$. We get $\#\Ad(\B)$ from~\cite[Theorem 5.12]{Ri1} and $\#\Ab(\B)$ from~\cite[Proposition 5.20(i),(ii),(v)]{Ri1}.

\vspace{0.25cm}
\noindent
$\G=\rD(m,n)$, $(m,n)\neq(2,1)$

\noindent
Note $\Delta=D_{m+2n}$, $r=s=n$, $S_j=\{\alpha_{i_j}\}$, $r_j=1$, $m_j=i_j$, and $l_j=i_j-j+1$ for $j=1,\ldots,n$. We get $\#\Ad(\B)$ from~\cite[Theorem 5.12]{Ri1} and $\#\Ab(\B)$ from~\cite[Proposition 5.20]{Ri1}.

\vspace{0.25cm}
\noindent
$\G=\rD(2,1;a)$

\noindent
Note $\Delta=D_4$, $r=s=1$, $S_1=\{\alpha_{i_1}\}$, $r_1=1$, $m_1=i_1$, and $l_1=i_1$. We get $\#\Ad(\B)$ from~\cite[Theorem 5.12]{Ri1} and $\#\Ab(\B)$ from~\cite[Proposition 5.20(i),(iv)]{Ri1}.

\vspace{0.25cm}
\noindent
$\G=\rF(4)$

\noindent
The result follows from the $E_7$ table in~\cite{Ri2}.

\vspace{0.25cm}
\noindent
$\G=\rG(3)$

\noindent
The result follows from the $E_8$ table in~\cite{Ri2}.
\end{proof}

Both (i) and (v) of~\cite[Proposition 5.20]{Ri1} reduce to $2^n+2^{n-1}-2$ in the $\rC(n)$ case. The ad-nilpotent equations for $\rC(n)$ also simplify to the closed forms given below. This is proven using Lemma~\ref{L:T_sum}, which we state after the proof of Proposition~\ref{P:C_closed}.

\begin{proposition}\label{P:C_closed}
Let $\G=\rC(n)$ and $n>2$. Then $\#\Ad(\B)$ is given by
\begin{enumerate}

\item[(i)] $(n+1)\cC_n-1$ if $c=0$,

\item[(ii)] $(4n-2)\cC_{n-1}+\binom{2n-1}{n-2}+\cT'_{n-1,n-2}-\cT'_{2n-c-2,c-1}-\cT'_{2n-c-2,c}$ if $c\in\{1,\ldots,n-1\}$.
\end{enumerate}
\end{proposition}

\begin{proof}
(i) It follows from Lemma~\ref{L:T_sum} that
\begin{equation}\label{e:C_sum1}
\sum_{j=1}^{n-1}\cT'_{2n-j-2,j}=\frac{n+1}{2}\cC_n-1.
\end{equation}
The result follows from the $c=0$ case of Table~\ref{Ta:enum_nonexcep}.

\noindent
(ii) It follows from Lemma~\ref{L:T_sum} that
\begin{equation}\label{e:C_sum2}
\sum_{j=1}^{n-1}\cT'_{2n-j-2,j-1}=\cT'_{n-1,n-2}+\binom{2n-1}{n-2}-\frac{n+1}{2}\cC_{n}+n\cC_{n-1}.
\end{equation}
Let $\cT_j=\cT'_{2n-j-2,j-1}+\cT'_{2n-j-2,j}$ for $j=1,\ldots,n-1$. Substituting~\eqref{e:C_sum1} and~\eqref{e:C_sum2}, we get
\begin{align*}
\sum_{j=c+1}^{n-1}\cT_j=\sum_{j=1}^{n-1}\cT_j-\sum_{j=1}^{c}\cT_{j}=\cT'_{n-1,n-2}+\binom{2n-1}{n-2}+n\cC_{n-1}-1-\sum_{j=2}^{c+1}\cT_{j-1}.
\end{align*}
The result follows from the $c\in\{1,\ldots,n-1\}$ case of Table~\ref{Ta:enum_nonexcep}.
\end{proof}

\begin{lemma}\label{L:T_sum}
Let $a,b,n\in\Z_{\geq 0}$ with $n\geq \lceil(a+3b)/2\rceil$ and $N=2n-a-b+2$. Then
\begin{equation*}
\sum_{j=0}^{\lfloor N/2\rfloor-1}\cT'_{2n-j-a,j-b}=\sum_{j=\lfloor N/2\rfloor-b-1}^{\lceil N/2\rceil+b}\binom{N-1} {j}-\frac{1}{2}\sum_{j=\lfloor N/2\rfloor-b}^{\lceil N/2\rceil+b}\binom{N}{j}.
\end{equation*}
\end{lemma}

\begin{proof}
We have the inequality $2n-j-a\geq j-b$ for $j=1,\ldots,\lfloor N/2\rfloor-1$, as
\begin{align*}
2n-\left\lfloor \frac{N}{2}\right\rfloor+1-a\geq \left\lfloor \frac{N}{2}\right\rfloor-1-b&&\text{if and only if}&& \left\lceil \frac{a+b}{2}\right\rceil\geq \frac{a-b}{2}.
\end{align*}
Furthermore, $n\geq \lceil(a+3b)/2\rceil$ implies that $\lfloor N/2\rfloor-1-b\geq 0$, and hence that $N\geq 1$ and $2n-j-a\geq 0$ for $j=1,\ldots, \lfloor N/2\rfloor-1$. Thus
\begin{equation*}
\sum_{j=0}^{\lfloor N/2\rfloor-1}\cT'_{2n-j-a,j-b}=\sum_{j=-b}^{\lfloor N/2\rfloor-1-b}\frac{N-2j}{N}\binom{N}{j}=\sum_{j=0}^{\lfloor N/2\rfloor-1-b}\binom{N}{j}-\frac{2}{N}\sum_{j=0}^{\lfloor N/2\rfloor-1-b}j\binom{N }{j}.
\end{equation*}
We simplify both summations in turn. Firstly,
\begin{align*}
\sum_{j=0}^{\lfloor N/2\rfloor-1-b}\binom{N }{j}&=\frac{1}{2}\left[\sum_{j=0}^{N}\binom{N}{j}-\sum_{j=\lfloor N/2\rfloor-b}^{\lceil N/2\rceil+b}\binom{N}{j}\right]=2^{N-1}-\frac{1}{2}\sum_{j=\lfloor N/2\rfloor-b}^{\lceil N/2\rceil+b}\binom{N}{j}.
\end{align*}
Secondly, and after recalling the $j\binom{N}{j}=N\binom{N-1}{j-1}$ identity,
\begin{equation*}
\frac{2}{N}\sum_{j=0}^{\lfloor N/2\rfloor-1-b}j\binom{N}{j}=\sum_{j=0}^{N-1}\binom{N-1}{j}-\sum_{j=\lfloor N/2\rfloor-1-b}^{\lceil N/2\rceil+b}\binom{N-1}{j}=2^{N-1}-\sum_{j=\lfloor N/2\rfloor-1-b}^{\lceil N/2\rceil+b}\binom{N-1}{j}.
\end{equation*}
The result follows.
\end{proof}

\end{document}